\documentclass[12pt]{article}

\usepackage{amsmath, amsthm, amsfonts, amssymb}
\theoremstyle{definition}

\numberwithin{equation}{section}
\theoremstyle{remark}

\newcommand{\lra}{\longrightarrow}

\title{{\bf Criteria of Spectral Gap  for Markov Operators}
\footnote{supported in part by
NNSFC(11131003), SRFDP, and the Fundamental Research Funds for the Central Universities.}}
\author{{\bf Feng-Yu Wang}\\
\footnotesize{School of Mathematical Sciences,
Beijing Normal
University, Beijing 100875, China}\\
 \footnotesize{Department of Mathematics,
Swansea University, Singleton Park, SA2 8PP, United Kingdom}\\ \footnotesize{wangfy@bnu.edu.cn, F.-Y.Wang@swansea.ac.uk}}
\begin{document}
\def\R{\mathbb R} \def\Z{\mathbb Z} \def\ff{\frac} \def\ss{\sqrt}
\def\dd{\delta} \def\DD{\Delta} \def\vv{\varepsilon} \def\rr{\rho}
\def\<{\langle} \def\>{\rangle} \def\GG{\Gamma} \def\gg{\gamma}
\def\lll{\lambda} \def\LL{\Lambda} \def\nn{\nabla} \def\pp{\partial}
\def\d{\text{\rm{d}}} \def\bb{\beta} \def\aa{\alpha} \def\D{\mathcal D}
\def\F{\mathcal F} \def\E{\mathcal E} \def\si{\sigma} \def\ess{\text{\rm{ess}}}
\def\beg{\begin} \def\beq{\begin{equation}} \def\neq{\end{equation}}
\def\Ric{\text{\rm{Ric}}} \def\Hess{\text{\rm{Hess}}}
\def\i{\text{\rm{i}}} \def\ii{\text{\rm{ii}}} \def\iii{\text{\rm{iii}}}
\def\iv{\text{\rm{iv}}} \def\v{\text{\rm{v}}} \def\vi{\text{\rm{vi}}}
\def\vii{\text{\rm{vii}}} \def\viii{\text{\rm{viii}}} \def\e{\text{\rm{e}}}
\def\ra{\Rightarrow} \def\lra{\Leftrightarrow} \def\C{\mathcal C}
\def\tt{\tilde} \def\P{\mathbb P} \def\ll{\lambda} \def\N{\mathbb N}
\def\oo{\omega} \def\OO{\Omega} \def\H{\mathcal H} \def\mN{\mu^{\N}}
\def\L{\mathbb L} \def\H{\mathcal H} \def\tt{\tilde} \def\gap{\text{\rm{gap}}}
\def\kk{\kappa}\def\B{\mathcal B}\def\EE{\mathbb E}

\maketitle \abstract{Let $(E,\mathcal F,\mu)$ be a probability
space, and let $P$ be a   Markov  operator on
$L^2(\mu)$ with $1$  a simple eigenvalue such that $\mu P=\mu$ (i.e. $\mu$ is an invariant probability measure of $P$). Then $\hat P:=\ff 1 2 (P+P^*)$ has a spectral gap, i.e.  $1$ is isolated in the spectrum of $\hat P$, if and only if
$$\|P\|_\tau:=\lim_{R\to\infty}  \sup_{\mu(f^2)\le 1}\mu\big(f(Pf-R)^+\big)<1.$$  This   strengthens a conjecture of Simon and H$\phi$egh-Krohn  on the spectral gap  for hyperbounded operators  solved   recently by L. Miclo in \cite{M}.
Consequently, for a symmetric, conservative, irreducible Dirichlet form on $L^2(\mu)$,    a Poincar\'e/log-Sobolev type inequality holds if and only if so does  the corresponding defective inequality.  Extensions to sub-Markov operators and non-conservative Dirichlet forms are   also presented. }

\date{}
\
\vskip 12pt\noindent
AMS 1991 Subject classification: 47D07, 60H10

\vskip 12pt\noindent Key Words: Spectral gap,  ergodicity, tail norm, Poincar\'e inequality.
 \vskip 2cm

\section{Introduction}

Let $(E,\F, \mu)$ be a probability space. Let
  $P$ be a   Markov operator on $L^2(\mu)$ (i.e. $P$ is a linear operator on $L^2(\mu)$ such that $P1=1$ and $f\ge 0$ implies $Pf\ge 0$) such that $\mu P=\mu$ (i.e. $\mu$ is an invariant probability measure of $P$).  Let $\hat P=\ff 1 2(P+P^*)$ be the additive symmetrization of $P$, where $P^*$ is the adjoint operator of $P$ on $L^2(\mu)$. Assuming that $1$ is a simple eigenvalue of $P$, we aim to investigate the existence of spectral gap of $\hat P$ (i.e. $1$ is an isolated point in $\si(\hat P)$, the spectrum of $\hat P$).

It is well known by the ergodic theorem that $1$ is a simple eigenvalue of $P$ if and only if $P$ is  ergodic, i.e.  for $f\in L^2(\mu)$,
$$\lim_{n\to\infty}\ff 1 n\sum_{k=1}^n P^kf= \mu(f)\ \ \text{holds\ in\ }L^2(\mu).$$

Moreover, the ergodicity (i.e. $1$ is a simple eigenvalue)  is also equivalent to  the $\mu$-essential irreducibility (or resolvent-positive-improving property, see \cite{Wu}) of $P$:
\beq\label{A0}\sum_{n=1}^\infty \mu(1_AP^n1_B)>0,\ \ \mu(A),\mu(B)>0.\end{equation}
Indeed, if for some $f\in L^2(\mu)$   with $\mu(f)=0$ and $\mu(f^2)=1$ such that $Pf=f$, then by the Jensen inequality we have $Pf^+\ge (Pf)^+=f^+$. But $\mu(P f^+)=\mu(f^+)$, we conclude that $Pf^+=f^+$. Then for any $\vv>0$ such that $\mu(f<-\vv),\mu(f>\vv)>0$, we have
$$\mu(1_{\{f<-\vv\}}P^n 1_{\{f>\vv\}})\le\ff 1{\vv^2} \mu(f^-P^nf^+)=\ff 1{\vv^2}\mu(f^-f^+)=0,\ \ n\ge 1,$$ so that \eqref{A0} does not hold. On the other hand, if there exists $A,B\in\F$ with $\mu(A),\mu(B)>0$ such that $\mu(1_BP^n 1_A)=0$ for all $n\ge 1$, then the class
$$\C:= \{0\le f\le 1: \ \mu(1_BP^nf)=0,\ n\ge 1\}$$ contains the non-trivial function $1_A$. Since the family is bounded in $L^2(\mu)$, we may take a  sequence $\{f_n\}\subset \C$ which   converges weakly to some $f\in \C$ such that    $\mu(f)=\sup_{g\in\C} \mu(g)$.  As $f\lor (Pf)$ is also in $\C$, we have $\mu((Pf)\lor f)=\mu(f)=\mu(Pf)$. Thus, $Pf=f$. Since
$\mu(f)\ge\mu(A)>0$ and $\mu(1_BPf)=0$, we conclude that $f-\mu(f)$ is a non-trivial eigenfunction of $P$ with respect to $1$, so that $1$ is not a simple eigenvalue of $P$.

When $P$ is symmetric,   the spectrum $\si(P)$ of $P$ is contained in $[-1,1].$ Then $P$ has a spectral gap if $P$ is ergodic and  $\si(P)\subset \{1\}\cup [-1,\theta]$ for some $\theta\in [-1,1);$ or equivalently, the Poincar\'e inequality
\beq\label{P0}\mu(f^2)\le \ff 1 {1-\theta}\mu(f(f-Pf))+\mu(f)^2,\  \ \ f\in L^2(\mu)\end{equation} holds.
When $P$ is non-symmetric, \eqref{P0} is equivalent to   $\si(\hat P|_{\H_0})\subset [-1,\theta]$, where   $\H_0:=\{f\in L^2(\mu):\ \mu(f)=0\}$. Thus, \eqref{P0} holds for some $\theta\in [-1,1)$ if and only if $\hat P$ has a spectral gap.

Recall that $P$ is called hyperbounded if for some $p>2$
$$\|P\|_{2\to p}:= \sup_{\mu(f^2)\le 1} \mu(|Pf|^p)^{\ff 1p}<\infty.$$
It was conjectured   by Simon and H$\phi$egh-Krohn \cite{SH}  that if $P$ is symmetric, ergodic and hyperbounded,  then it has a spectral gap. Although  numerous papers aiming to solve this problem or to construct counterexamples have been published, see e.g. \cite{A,C,H,W04,Wu, GW} where some weaker notions such as the uniform integrability and a tail norm condition have been used to replace the hyperboundedness, the conjecture has been open for more than 40 years until Miclo found a complete proof in his recent paper \cite{M}.

On the other hand, there are a lot of non-hyperbounded Markov operators having   spectral gap. So, in the spirit of \cite{H,GW}, we shall prove a stronger statement by using a tail norm condition to replace the hyperboundedness. The tail norm  we will use is the following:
$$\|P\|_\tau:=\lim_{R\to\infty} \sup_{ \mu(f^2)\le 1} \mu\big(f(Pf-R)^+\big)= \lim_{R\to\infty} \sup_{f\ge 0, \mu(f^2)\le 1} \mu\big(f(Pf-R)^+\big).$$
According to  the Schwartz inequality, $\|P\|_\tau$ is smaller  than the following one used in \cite{H,GW}:
$$\|P\|_{tail}:= \lim_{R\to\infty} \sup_{ \mu(f^2)\le 1} \mu\big({(|Pf|-R)^+}^2\big)^{\ff 1 2}.$$  Moreover,  by the Jensen inequality
$\|P^m\|_{tail}$ is decreasing in $m\in\N$. When $P$ is symmetric,  $\|P^{2m-1}\|_{\tau}$ has the same property     since for any $R>0$ and $m\in\N$,
 \beg{equation*}\beg{split}& \sup_{\mu(f^2)\le 1} \mu(f(P^{2m+1}f-R)^+) \le \sup_{\mu(f^2)\le 1} \mu(fP(P^{2m}f-R)^+)\\
 &=\sup_{\mu(f^2)\le 1} \mu((Pf)(P^{2m-1}Pf-R)^+)
\le \sup_{\mu(f^2)\le 1} \mu(f(P^{2m-1}f-R)^+),\end{split}\end{equation*}where the las step follows by replacing $Pf$ with $f$, since $\mu(f^2)\le 1$ implies $\mu(Pf)^2)\le 1$.

According to \cite{Wu},  $P$ is called  uniformly integrable  if $\|P\|_{tail}=0$. Thus, the uniformly integrability implies  $\|P\|_\tau=0.$ In particular,
$\|P\|_\tau=0$ holds for hyperbounded $P$. We will then strengthen the above conjecture by replacing the hyperboundedness with $\|P\|_\tau<1,$ which is also necessary for the existence of the spectral gap of $\hat P$ as shown in our following main result. See Theorem \ref{T2.1} below  for two more equivalent statements on isoperimetric constants for  the existence of spectral gap.

\beg{thm}\label{T1.1} Let $P$ be an ergodic  Markov operator on $L^2(\mu)$. Then the following statements are equivalent:
\beg{enumerate} \item[$(1)$] The Poincar\'e inequality $\eqref{P0}$ holds for some constant $\theta\in [-1,1)$, i.e. $\hat P$ has a spectral gap.
\item[$(2)$]   $\|P\|_\tau<1.$
\item[$(3)$]   $\inf_{m\in\N}\|P^{2m-1}\|_\tau<1$.
\item[$(4)$]   $\inf_{m\in\N}\|{\hat P}^m\|_{tail}<1.$
\ \end{enumerate}\end{thm}

\paragraph{Remark 1.1.} (a) Let $P$ be symmetric. It is easy to see that if $P^2$ has a spectral gap then $\|P\|_{tail}<1$. Indeed, letting $\si(P^2)\subset [0,\vv]$ for some $\vv\in (0,1)$, for $R\ge 1$ and $\mu(f^2)\le 1$ we have
$$\mu({(Pf-R)^+}^2)\le \mu(|P\hat f|^2)= \mu(\hat fP^2\hat f)\le\vv\mu({\hat f}^2)\le\vv,$$ where $\hat f:= f-\mu(f).$ On the other hand, if $\|P\|_{tail}<1$ then by Theorem \ref{T1.1} $P^2$ has a spectral gap provided it is ergodic (i.e. $1$ is a simple eigenvalue). So, if $P$ has a spectral gap but $-1
\in\si_{ess}(P)$ is not an eigenvalue, then $\|P\|_{tail}=1$. From this we would believe that in general   the existence of spectral gap for $P$ does not imply $\|P\|_{tail}<1.$



(b) The equivalence of the existence of spectral gap   and $\inf_{m\in\N}\|P^m\|_{tail}<1$  has been proved in \cite{GW} for resolvent-uniform-positive-improving Markov operators (see \cite[Lemma 3.6,Theorem 4.1]{GW}): for any $\vv>0$, there exists $m\in \N$ such that
$$\inf\Big\{\sum_{k=1}^m\mu(1_AP^k1_B): \ \mu(A),\mu(B)\ge\vv\Big\}>0.$$  The resolvent-uniform-positive-improving condition is weaker than the uniform-positive-improving condition used in \cite{A,K}, and the  condition (E) introduced in \cite{H}. Recall that $P$ is called uniform-positive-improving if
$$\inf\{\mu(1_AP1_B):\ \mu(A),\mu(B)\ge\vv\}>0,\ \ \vv\in (0,1].$$

In the symmetric setting this result is now improved by Theorem \ref{T1.1} since   the resolvent-uniform-positive-improving property  is strictly stronger than the ergodicity, which is equivalent to the   resolvent-positivity-improving property (\ref{A0}). Examples to ensure this ``strictly stronger" property can be constructed from the corresponding ones for Markov semigroups explained in item (c) below. More precisely, it is known that the resolvent-uniform-positive-improving property of $P$ is equivalent to the validity of the weak Poincar\'e inequality (see \cite[Theorem 2.5]{Wu2}):
 $$\mu(f^2)\le\aa(r) \mu(f(1-P)f)+r\|f\|_\infty^2,\ \ r>0, f\in L^2(\mu), \mu(f)=0 $$ for some $\aa: (0,\infty)\to (0,\infty).$ Taking $P=P_1$ for a symmetric Markov semigroup $P_t$ on $L^2(\mu)$ and noting that
 $$\mu(f(1-P_1)f)=\mu(f(1-\e^{L})f)\le -\mu(fLf)=\E(f,f),\ \ f\in\D(L),$$ where $(L,\D(L))$ and $\E$ are the associated generator and Dirichlet form, we see that an example such that $\E$ is irreducible but the weak Poincar\'e inequality \eqref{AD} below is not available implies that $P$ is ergodic but does not possess the resolvent-uniform-positive-improving property.

(c) We would like to mention links of the uniform-positive-improving property of  a symmetric Markov semigroup $P_t$ and the weak Poincar\'e inequality of the associated Dirichlet form $(\E,\D(\E))$. It is well known that $P_t$ (for some/all $t>0$) is ergodic if and only if the Dirichlet form is irreducible, i.e. $\E(f,f)=0$ implies $f$ is constant. Next,  according to \cite{K} (see also \cite{A}), the uniform-positive-improving property of the semigroup (i.e. $P_t$ is uniform-positive-improving for some, equivalently all,  $t>0,$ see \cite{H}) implies  the weak spectral gap property, which is equivalent to the validity of the weak Poincar\'e inequality (see \cite{RW}): for some $\aa:(0,\infty)\to (0,\infty)$
\beq\label{AD} \mu(f^2)\le \aa(r)\E(f,f)+ r\|f\|_\infty^2,\ \ r>0, f\in\D(\E), \mu(f)=0.\end{equation}  Indeed, as explained in (b) that the resolvent-uniform-positive-improving property of $P_t$ for some $t>0$ also implies the existence of the weak Poincar\'e inequality. Moreover, it is shown in  \cite[\S7]{RW} that  there are conservative irreducible Dirichlet forms  which do not satisfy the weak Poincar\'e inequality.  Therefore, the resolvent-uniform-positive-improving property is strictly stronger than the ergodicity.



\

As applications of Theorem \ref{T1.1}, we consider functional inequalities conservative  symmetric Dirichlet forms. A simple consequence of the equivalence of $(1)$ and (3) is that the defective Poincar\'e inequality implies the tight one.

\beg{cor}\label{C1.2} Let $(\E,\D(\E))$ be a conservative, irreducible, symmetric Dirichlet form on $L^2(\mu)$. Then the Poincar\'e inequality
\beq\label{P} \mu(f^2)\le C\E(f,f)+ \mu(f)^2,\ \ \ f\in\D(\E)\end{equation}  holds for some constant $C>0$ if and only if the defective Poincar\'e inequality
\beq\label{DP} \mu(f^2)\le C_1\E(f,f)+C_2\mu(|f|)^2,\ \ f\in\D(\E)\end{equation} holds for some constants $C_1,C_2>0.$
\end{cor}
This result improves   \cite[Proposition 1.3]{RW} where the weak Poincar\'e inequality (\ref{AD})  is used to replace the irreducibility of the Dirichlet form.  Basing on Corollary \ref{C1.2}, we are able to prove the equivalence of the defective version and the tight version for more general functional inequalities. Here, we consider a family of functional inequalities introduced in \cite{LO}, which interpolate the Poincar\'e inequality (\ref{P}) and the Gross \cite{Gross} log-Sobolev inequality
\beq\label{LS} \mu(f^2\log f^2)\le C\E(f,f),\ \ f\in\D(\E), \mu(f^2)=1\end{equation} for some constant $C>0.$  Let $\phi\in C([1,2])$ such that $\phi>0$ on $[1,2)$ and $\phi(2)=0$ with
$$c_\phi:= \sup_{p\in [1,2)} \ff{2-p}{\phi(p)}<\infty.$$ Consider the functional inequality
\beq\label{PH} {\rm Var}_{\mu,\phi}(f):= \sup_{p\in [1,2)} \ff{\mu(f^2)-\mu(|f|^p)^{\ff 2 p}}{\phi(p)} \le C\E(f,f),\ \ f\in\D(\E)\end{equation} and its defective version
\beq\label{DPH} {\rm Var}_{\mu,\phi}(f)\le C_1\E(f,f)+C_2\mu(f^2),\ \ f\in\D(\E).\end{equation} When $\phi(p)=2-p$ the inequality (\ref{PH}) is equivalent to the log-Sobolev inequality, and when $\phi$ reduces to a positive constant it becomes the Poincar\'e inequality. See \cite{W05} for detailed discussions on properties and applications of the inequality (\ref{PH}).

\beg{cor}\label{C1.3} Let $(\E,\D(\E))$ be a conservative, irreducible, symmetric Dirichlet form on $L^2(\mu)$. Then $(\ref{PH})$ holds for some constant $C>0$ if and only if $(\ref{DPH})$ holds for some constants $C_1,C_2>0.$\end{cor}

The remainder of the paper is organized as follows. In Section 2, by using an approximation argument introduced in \cite{M},  we extend
a known Cheeger type inequality  for high order eigenvalues of finite-state Markov chains to the abstract setting, then use this estimate to characterize the existence of spectral gap with high-order isoperimetric constants. This characterization is then used in Section 3  to prove Theorem \ref{T1.1}. Proofs of        Corollaries \ref{C1.2} and \ref{C1.3}  are also addressed in Section 3. Finally, in Section 4 we extend Theorem \ref{T1.1} and Corollary \ref{C1.2} to the sub-Markov setting.

\section{Essential spectrum   and isoperimetric constants}
Let $P$ be a symmetric Markov operator on $L^2(\mu)$. We aim to characterize the essential spectrum of $P$ using high order isoperimetric constants
investigated in \cite{24}.

For any $n\in\N$,  let
$$D_n=\{(A_1,\cdots,A_n):\ A_1,\cdots A_n\ \text{are\ disjoint},\ \mu(A_k)>0,1\le k\le n\}.$$We define the $n$-th isoperimetric constant by
$$\kk_n= \inf_{(A_1,\cdots,A_n)\in D_n} \max_{1\le k\le n} \ff{\mu(1_{A_k}P1_{A_k^c})}{\mu(A_k)}.$$
Obviously, $\kk_1=0$ and   if $L^2(\mu)$ is infinite-dimensional then $D_n\ne\emptyset$ for all $n\ge 2.$
It is also easy to see that $\kk_n$ is non-decreasing in $n$.

We will only consider the case that  $L^2(\mu)$ is infinite-dimensional, since otherwise   the spectrum of $P$ is finite so that the existence of spectral gap becomes  trivial.    By the Cheeger inequality we know that $P$ has a spectral gap if and only if $\kk_2>0$, see e.g. \cite{LS} by noting that in (\ref{P0}) we have
$$\mu(f(f-Pf))=\ff 1 2 \int_{E\times E} (f(x)-f(y))^2J(\d x,\d y)$$   for the symmetric measure $J$ on $E\times E$ determined by
$J(A\times B):= \mu(1_AP1_B),A,B\in\F$.

Consider $\ll_{ess}(P):=\sup\si_{ess}(P),$ where $\si_{ess}(P)$ is the essential spectrum of $P$.
Obviously, $P$ has a spectral gap if and only if it is ergodic and $\ll_{ess}(P)<1.$

\beg{thm}\label{T2.1} Let $P$ be a symmetric  Markov operator on $L^2(\mu)$. Then   $\ll_{ess}(P)<1$ if and only if
 $\sup_{n\ge 1} \kk_n>0.$  Consequently, when $P$ is ergodic then the existence of spectral gap of $P$ is equivalent to each of the following two statements:\beg{enumerate}
  \item[$(5)$] $\ \kk_2 \big(=\inf_{n\ge 2}\kk_n\big)>0.$  \item[$(6)$] $ \lim_{n\to\infty}\kk_n \big(=\sup_{n\ge 2}\kk_n\big)>0.$\end{enumerate}\end{thm}

As mentioned at the end of Introduction, to prove this result we will extend  a known estimate on the hight order eigenvalues using $\kk_n$ for finite-state Markov chains. So, below we first consider Markov operators on a finite set.

Let $n\ge 1$ be fixed, and let $\tt E$ be a finite set with $|\tt E|\ge n$, where $|\tt E|$ denotes the number of elements in $\tt E$. Let $\tt\mu$ be a strictly positive probability measure on $\tt E$ equipped with the largest $\si$-field $\B(\tt E)$, i.e. $\B(\tt E)$ is the class of all subsets of $\tt E$ and $\tt\mu(\{x\})>0$ for any
$x\in\tt E$. For a symmetric Markov operator $\tt P$ on $\tt E$, let
$$0=\tt \ll_1\le\cdots\le \tt\ll_{|\tt E|}\le 2$$ be all eigenvalues of $1-\tt P$. According to \cite{24}, there exists a constant $c(n)>0$  depends only on $n$ such that
\beq\label{B1} \ss{\tt\ll_n}\ge c(n)\inf_{(\tt A_1,\cdots, \tt A_n)\in\tt D_n} \max_{1\le k\le n} \ff{\tt\mu(1_{\tt A_k} \tt P 1_{\tt A_k^c})}{\tt\mu(\tt A_k)},\end{equation} where
$$\tt D_n=\{(\tt A_1,\cdots, \tt A_n): \tt A_1,\cdots \tt A_n \ \text{are\ disjoint\ non-empty\ subsets\ of\ } \tt E\}.$$
As shown in \cite[Theorem 2]{M}, we may take $c(n)= \ff{c_0}{n^4}$ for a universal constant $c_0>0.$

Below we aim to extend this estimate to our abstract setting by using an approximation argument introduced  in \cite{M}. For any $n\ge 1$, let
\beq\label{**}\ll_n= \sup_{f_1,\cdots, f_{n-1}\in L^2(\mu)} \inf_{{\mu(f^2)=1,}\atop{\mu(ff_i)=0, 1\le i\le n-1}} \mu\big(f(1-P)f\big).\end{equation}
To see that  this quantity can be regarded as the $n$-th eigenvalue of $L:=1-P$, let $\ll_{ess}(L)=\inf\si_{ess}(L)$ be the bottom  of the essential spectrum of $L$. Then, see e.g.  \cite[Theorem XIII.2]{RS}, $\ll_n$ is the $n$-th eigenvalue of $L$ if $\ll_n<\ll_{ess}(L),$ and $\ll_n=\ll_{ess}(L)$ otherwise. The following result was stated as Proposition 5 in \cite{M}, we include here a   proof for completeness.

\beg{lem}[\cite{M}]\label{L2.2} Let $c(n)$ be in $(\ref{B1})$ for $n\ge 1$. Then $\kk_n\ge \ll_n\ge c(n)^2\kk_n^2.$\end{lem}

\beg{proof}  Let $L=1-P$. Since $\ll_1=\kk_1=0$, we only prove for $n\ge 2.$

(a) Upper bound estimate of $\ll_n$. For any $f_1,\cdots, f_{n-1}\in L^2(\mu)$ and any $(A_1,\cdots, A_n)\in D_n$, there exists a function
$f:= \sum_{i=1}^n a_i 1_{A_i}$ such that $\sum_{i=1}^n a_i^2=1$ and $\mu(ff_i)=0, 1\le i\le n.$ Let $$\kk= \max_{1\le i\le n} \ff{\mu(1_{A-i}P1_{A_i^c})}{\mu(A-i)}.$$Then
\beg{equation*}\beg{split} &\mu(fLf)= \sum_{i,j=1}^n a_ia_j \mu(1_{A_i}L1_{A_j})\le \sum_{i,j=1}^n |a_i|\cdot |a_j|\ss{\mu(1_{A_i}L1_{A_i})\mu(1_{A_j}L1_{A_j})}\\
& \le \kk\sum_{i,j=1}^n |a_i|\cdot |a_j|\ss{\mu(A_i)\mu(A_j)}
\le \kk\sum_{i,j=1}^na_i^2\mu(A_j)\le\kk.\end{split}\end{equation*}
Therefore, by the definition of $\ll_n$ and $\kk_n$, we have $\ll_n\le\kk_n$.

(b) Lower bound estimate of $\ll_n.$ Assume that $\ll_n<c(n)^2\kk_n^2$. Then there exist
$f_1,\cdots,f_{n}\in L^2(\mu)$ such that
\beq\label{A*0} \beg{split} &\mu(f_if_j)=\dd_{ij},\ \ \mu(f_iLf_i)\le c(n)^2\kk_n^2-\ff{\dd_n}2,\\
&\max_{1\le i\le n} \sum_{j\ne i,1\le j\le n}|\mu(f_iLf_j)|\le \ff {\dd_n} 4,\ \ \le i,j\le n,\end{split}\end{equation}
where $\dd_n:= c(n)^2\kk_n^2-\ll_n>0.$  Indeed, for any $\vv>0$ we may find an orthonormal family $\{f_1,\cdots,f_n\}\subset L^2(\mu)$ such that $\mu(|Lf_i-\ll_i f_i|^2)\le\vv, 1\le i\le n.$ Since $\max_{1\le i\le n}\ll_i=\ll_n=c(n)^2\kk_n^2-\dd_n$ and $\mu(f_if_j)=\dd_{ij}$,  (\ref{A*0}) holds for small enough $\vv>0.$

Let $\F_\infty=\si(f_1,\cdots, f_n).$ Since $\F_\infty$ is separable, we may find an increasing sequence of $\si$-fields $\{\F_N\}_{N\ge 1}$ such that
 $$\bigvee_{N\ge 1} \F_N= \F_\infty.$$
 Let $\mu_N$ be the restriction of $\mu$ on $\F_N$, and let $\EE^{\F_N}$ be the conditional expectation under $\mu$ given $\F_N$. Let
 $$P_N g= \EE^{\F_N} Pg,\ \ g\in L^2(\mu_N).$$ Then $P_N$ is a symmetric Markov operator on $L^2(\mu_N).$

  To identify $P_N$ with a Markov operator on a finite set, let
 $$\tt E_N=\{A: A\text{\ is\ atom\ of\ }\F_N, \mu_N(A)>0\},\ \ \tt\mu_N(\{A\})=\mu_N(A)=\mu(A),\ \ A\in\tt E_N.$$ Then $\tt E_N$ is a finite set with $\tt\mu_N$
 a strictly probability measure. It is easy to see that the map
 $$\varphi_N: \tt g\mapsto \sum_{A\in \tt E_N} \tt g(A) 1_A$$ is   isometric from $L^2(\tt\mu_N)$ to $L^2(\mu_N).$ Moreover, the inverse of $\varphi_N$ is given by (note that $g$ is constant on atoms of $\F_N$)
 $$\varphi^{-1}_N(g)(A)= g|_A,\ \ A\in \tt E_N,\ g\in L^2(\mu_N).$$ Define
 $$\tt P_N\tt g=\varphi_N^{-1}(P_N\varphi_N(\tt g)),\ \ \tt g\in L^2(\tt\mu_N).$$ Then $\tt P_N$ is a symmetric Markov operator on $L^2(\tt\mu_N)$ and having the same spectral information of $P_N$. Therefore, (\ref{B1}) is valid for $P_N$, i.e. letting
 $$0=\ll_{N,1}\le\cdots\le \ll_{N,\ell_N}\le 2$$ be all eigenvalues of $L_N:=1-P_N$, where $\ell_N=|\tt E_N|={\rm dim}L^2(\mu_N),$
 if $\ell_N\ge n$ then
 $$ \ss{\ll_{N,n}}\ge c(n)\inf_{(A_1,\cdots, A_n)\in D_{N,n}} \max_{1\le k\le n} \ff{\mu_N(1_{ A_k}  P_N 1_{ A_k^c})}{\mu_N(A_k)}$$ holds for
 $D_{N,n}:= \{(A_1,\cdots, A_n)\in D_n:\ A_k\in\F_N, 1\le k\le n\}.$ Since for $A_k\in\F_N$ we have $\mu_N(A_k)=\mu(A_k)$ and
 $$\mu_N(1_{A_k}P_N1_{A_k^c})= \mu(1_{A_k} \EE^{\F_N}P1_{A_k^c})=\mu(1_{A_k}P1_{A_k^c}),$$ this implies
 \beq\label{B3}  \ll_{N,n}\ge c(n)^2\kk_n^2,\ \ \  n\le \ell_N.\end{equation}

  Now, let $f_{N,k}= \EE^{\F_N} f_k, 1\le k\le n.$ Then by the martingale convergence theorem $\lim_{N\to\infty} f_{N,k}=f_k$ holds in $L^2(\mu), 1\le k\le n.$ Let $L_N=1-P_N$. Since
 $P$ is continuous in $L^2(\mu)$, combining this with  (\ref{A*0}) we obtain
\beg{equation*}\beg{split} & \lim_{N\to\infty} \mu_N(f_{N,i}f_{N,j})=\lim_{N\to\infty} \mu(f_{N,i}f_{N,j})= \dd_{ij},\ \ 1\le i,j\le n,\\
 &\lim_{N\to\infty} \mu(f_{N,i}L_Nf_{N,i})  =\lim_{N\to\infty}  \mu(f_{N,i}(1-P)f_{N,i})\le c(n)^2\kk_n^2-\ff{\dd_n}2,\ \ 1\le i\le n,\\
 &\lim_{N\to\infty} \sup_{1\le i\le j} \sum_{j\ne i}    |\mu(f_{N,i}L_Nf_{N,j})|\le \ff{\dd_n}4.\end{split}\end{equation*} This is contradictive to (\ref{B3}). Indeed, from this we may find orthonormal $\{\tt f_{N,i}: 1\le i\le n\}\subset L^2(\mu_N)$ such that
 $$\vv_N:=\max_{1\le i\le n} \mu(|f_{N,i}-\tt f_{N,i}|^2)^{\ff 1 2 } \to  0\ \text{as}\ N\to\infty,$$ and thus,
 \beg{equation*}\beg{split}
 &\mu(\tt f_{N,i}L_N \tt f_{N,i})\le \mu(f_{N,i}L_N f_{N,i})+4\vv_N\le c(n)^2\kk_n^2-\ff{\dd_n} 2 +4\vv_N,\ \ 1\le i\le n,\\
 &\max_{1\le i\le n}\sum_{j\ne i} |\mu(\tt f_{N,i}L_N\tt f_{N,j})|\le\ff{\dd_n} 4 + 4(n-1)\vv_N.\end{split}\end{equation*} Therefore, for $f:=\sum_{i=1}^n a_i\tt f_{N,i}$ with $\sum_{i=1}^n a_i^2=1,$
 \beg{equation*}\label{**2} \beg{split} &\mu_N(fL_Nf)= \sum_{i,j=1}^n \mu(\tt f_{N,i}L_N\tt f_{N,j})a_ia_j\\
 &\le c(n)^2\kk_n^2-\ff{\dd_n}2+4\vv_N+\ff{\dd_n}4 +4(n-1)\vv_N= c(n)^2\kk_n^2 -\ff{\dd_n}4 +4n\vv_N.\end{split}\end{equation*}
 Since for any $g_1,\cdots,g_{n-1}\in L^2(\mu_N)$ there exists $f\in\text{span}\{\tt f_{N,i}: 1\le i\le n\}$ with $\mu(f^2)=1$ such that $\mu_N(fg_i)=0, 1\le i\le n-1,$ combining this with (\ref{**}) for $P_N$ in place of $P$, we arrive at
 $$\ll_{N,n}\le c(n)^2\kk_n^2 -\ff{\dd_n}4 +4n\vv_N,$$   which contradicts  (\ref{B3}) for large $N$ such that $\vv_N<\ff{\dd_n}{16n}$.  \end{proof}

  \beg{proof}[Proof of Theorem \ref{T2.1}] If $\ll_{ess}(P)<1$, then
    $\si(L)\cap [0, 1-\ll_{ess}(P))$ is discrete and each eigenvalue in this set is of finite multiplicity. So, in this case  $\ll_n>0$ for $n$ lager than the multiplicity of the first eigenvalue $\ll_1=0,$ and hence  by Lemma \ref{L2.2},
  $  \kk_n>0 $  for large $n$. On the other hand,   we aim to prove that if $\ll_{ess}(P)=1$  then $\kk_n=0$ for all $n\ge 1.$ Since $\ll_{ess}(P)=1$, i.e. $0\in\si_{ess}(L)$,    we  have $\ll_n= 0$ for all $n$. Combining this with   Lemma \ref{L2.2}, we prove $\kk_n=0$ for all $n\ge 1$. Thus, the proof of the first assertion is finished.

  Now, let $P$ be ergodic. Since (6) follows from (5), it suffices to prove that (1) (i.e. the existence of spectral gap of $P$) implies (5) while (6) implies (1). If (1) holds then $\ll_2>0$, so that by Lemma \ref{L2.2} we have $\kk_2>0$, i.e. (5) holds. On the other hand, if (6) holds, then by the first assertion of this theorem we have $\ll_{ess}(P)<1$ such that $1$ is isolate in $\si(P)$. Since $1$ is a simple eigenvalue of $P$, this implies that $P$ has a spectral gap, i.e. (1) holds.
  \end{proof}

\section{Proofs of Theorem \ref{T1.1} and corollaries }




\beg{proof}[Proof of Theorem \ref{T1.1}] By \cite[Theorem 4.1(d)]{GW} with $\pi=\hat P$ and $r_{sp}(\pi)=1$, (1) implies (4). Next, it is obvious that (2) implies (3). Moreover, since $\|\hat P^m\|_{\tau}\le \|\hat P^m\|_{tail}$ and by the Jensen inequality the later is decreasing in $m$, we see that   (4)   implies (3) for $\hat P$ in place of  $P$.    Therefore, it remains to prove that
  (1) implies (2), and (3) implies (1). Below we prove these two implications respectively.


 {\bf (1) implies (2). }  Let $\theta\in (0,1)$ such that (\ref{P0}) holds. For any $f\ge 0$ with $\mu(f^2)\le 1$, we have
$$\mu(f^2)\le \ff 1 {1-\theta}\mu(f(f-Pf)) +\mu(f)^2= \ff{\mu(f^2)}{1-\theta}-\ff{\mu(fPf)}{1-\theta} +\mu(f)^2.$$
This implies that
$$\mu(fPf)\le \theta\mu(f^2)+ (1-\theta)\mu(f)^2\le \theta +(1-\theta)\mu(f)^2.$$ Replacing $f$ by $(f-\ss R)^+$ for $R>1$, and noting that
\beq\label{DD} \mu((f-s)^+)^2\le \mu\big({(f-s)^+}^2\big)\mu(f>s)\le \ff 1 {s^2},\ \ s>0,\end{equation}
we obtain
$$\mu\Big(\big(f-\ss R\big)^+ P\big(f-\ss R\big)^+\Big)\le \theta +\ff{1-\theta}R.$$ Since by the Jensen inequality $(Pf-R)^+\le P(f-R)^+,$ combining this with
(\ref{DD}) we arrive at
\beg{equation*}\beg{split} &\mu\big(f(Pf-R)^+\big)\le \mu\big(fP(f-R)^+\big)\\
&\le \mu\Big(\big(f-\ss R\big)^+P(f-R)^+\Big) +\ss R\, \mu\big(P(f-R)^+\big)\\
&\le\mu\Big(\big(f-\ss R\big)^+ P\big(f-\ss R\big)^+\Big) +\ss R\, \mu\big((f-R)^+\big)\\
&\le \theta+ \ff{1-\theta}R+\ff 1 {\ss R},\ \ R>1.\end{split}\end{equation*}
Therefore,
$$\|P\|_\tau= \lim_{R\to\infty} \sup_{0\le f, \mu(f^2)\le 1} \mu\big(f(Pf-R)^+\big)\le\theta<1.$$

{\bf (3) implies (1)}. Let   $\|P^{2m-1}\|_\tau<1$ for some $m\ge 1$. It suffices to prove for the case that $L^2(\mu)$ is infinite-dimensional since in the finite dimensional case the existence of spectral gap is trivial. We first assume that $P$ is symmetric. In this case,   the existence of spectral gap for $P$ is equivalent to that for $P^{2m-1}$.   Since  $\|P^{2m-1}\|_\tau<1$, there exists $\dd\in (0,1)$ and $R>0$  such that
\beq\label{D1} \sup_{f\ge 0, \mu(f^2)\le 1} \mu(f(P^{2m-1}f-R)^+)\le \dd.\end{equation}
On the other hand, if $P^{2m-1}$ does not have spectral gap, by Theorem \ref{T2.1}, for any $\vv\in (0,1)$ and   $n\ge 2$, there exists $(A_1,\cdots A_n)\in D_n $
such that
$$\max_{1\le k\le n} \ff{\mu(1_{A_k}P^{2m-1}1_{A_k^c})}{\mu(A_k)}\le\vv.$$ Letting $B_k=A_k\cap \{P^{2m-1}1_{A_k}\ge 1-\ss\vv\}$, this implies
\beg{equation*}\beg{split} &\vv\mu(A_k)\ge \mu(1_{A_k}P^{2m-1}1_{A_k^c}) =\mu(A_k) -\mu(1_{A_k}P^{2m-1}1_{A_k})\\
&\ge \mu(A_k)-\mu(B_k)-\big(1-\ss\vv\big)\big(\mu(A_k)-\mu(B_k)\big).\end{split}\end{equation*} Thus,
\beq\label{D2} \mu(A_k)\ge\mu(B_k) \ge \big(1-\ss\vv\big)\mu(A_k).\end{equation}
Since $A_1,\cdots, A_n$ are disjoint with positive $\mu$-mass  and $\mu$ is a probability measure, there exists $1\le k\le n$ such that $\mu(A_k)\in (0,\ff 1 n).$ Take
$$f=\ff{1_{A_k}}{\ss{\mu(A_k)}}.$$ We have
$$P^{2m-1}f= \ff{P^{2m-1}1_{A_k}}{\ss{\mu(A_k)}}\ge \ff{1-\ss\vv}{\ss{\mu(A_k)}}1_{B_k}.$$ Combining this with (\ref{D1}), (\ref{D2}) and $\mu(A_k)\in (0,\ff 1 n)$, we arrive at
$$\dd\ge \mu(f(P^{2m-1}f-R)^+) \ge \ff{\mu(B_k)\big(1-\ss\vv-R\ss{\mu(A_k)}\big)^+}{\mu(A_k)}\ge \big(1-\ss\vv\big)\Big(1-\ss\vv-\ff R{\ss n}\Big)^+.$$
Since $\vv\in (0,1)$ and $n\ge 2$ are arbitrary, by letting $\vv\to 0$ and $n\to\infty$ we obtain $\dd\ge 1$ which  contradicts  $\dd\in (0,1)$.

Next, for non-symmetric $P$, we consider the symmetrizing operator $\hat P$. Obviously, $P$ satisfies (\ref{A0}) implies that   $\hat P$ satisfies (\ref{A0}), i.e. $\hat P$ is ergodic. By (3), there exists $m\ge 1$ such that   $\vv:=\|P^{2m-1}\|_\tau<1$. Moreover, as proved above if $\|{\hat P}^{2m-1}\|_\tau<1$ then $\hat P$ has a spectral gap. Therefore, it suffices to show that $\|{\hat P}^{2m-1}\|_\tau<1.$  Noting that for any $R>0$ and  $f\ge 0$ with $\mu(f^2)\le 1$,
\beg{equation*}\beg{split} \mu\big(f({\hat P}^{2m-1}f-R)^+ \big)&\le
2^{1-2m} \sum_{i_1,\cdots, i_{2m-1}=0}^1\mu\Big(f  \Big(\Big\{\prod_{k=1}^{2m-1} P^{i_k}(P^*)^{1-i_k}\Big\}f-R\Big)^+\Big)\\
&\le 2^{1-2m} \mu\big(f (P^{2m-1}f-R)^+\big)+ 2^{1-2m} \big(2^{2m-1}-1\big),\end{split} \end{equation*}
we obtain $\|{\hat P}^{2m-1}\|_\tau \le 2^{1-2m}\vv + 1- 2^{1-2m}<1$, since $\vv<1.$\end{proof}

\beg{proof}[Proof of Corollary \ref{C1.2}] It suffices to prove (\ref{P}) from (\ref{DP}). Let $P_t$ be the associated Markov semigroup. Then the irreducibility of the Dirichlet form implies that of $P_t$ for $t>0$. Next, by
\cite[Theorem 3.3]{W02} with $\phi\equiv 1$, (\ref{DP}) implies that $\|P_t\|_\tau\le \e^{-t/C_1}<1$ for $t>0$. Then due to Theorem \ref{T1.1} we conclude that $P_t$ has a spectral gap, equivalently, the Poincar\'e inequality (\ref{P}) holds for some constant $C>0.$\end{proof}

\beg{proof}[Proof of Corollary \ref{C1.3}]
 By \cite[Proposition 2.1]{W05}, for any $f\in L^2(\mu)$ and $\hat f:=f-\mu(f)$, we have
\beg{equation*}\beg{split} &\mu(f^2)-\mu(|f|^p)^{\ff 2 p}\le \mu({\hat f}^2) + (1-p) \mu(|\hat f|^p)^{\ff 2 p}\\
&= (2-p) \mu({\hat f}^2)+(p-1) \big(\mu({\hat f}^2)-\mu(|\hat f|^p)^{\ff 2 p}\big).\end{split}\end{equation*}
Then it follows from (\ref{DPH}) that
\beq\label{*0} {\rm Var}_{\mu,\phi}(f)\le c_\phi\mu({\hat f}^2) +{\rm Var}_{\mu,\phi}(\hat f) \le (c_\phi+C_2)\mu({\hat f}^2) + C_1\E(f,f).\end{equation} Next, since $\phi(2)=0$, from the proof of \cite[Theorem 1.1(2)]{W05} we see that the $F$-Sobolev inequality (1.7) in \cite{W05} holds for some nonnegative function $F$ with $F(r)\uparrow\infty$ as $r\uparrow\infty$. According to \cite{W00} (see also \cite{Wbook}), this inequality is equivalent to the super Poincar\'e inequality
$$\mu(f^2)\le r\E(f,f) +\bb(r) \mu(|f|)^2,\ \ r>0, f\in \D(\E)$$ for some function $\bb: (0,\infty)\to (0,\infty)$. In particular, the defective Poincar\'e inequality holds. Thus, by Corollary \ref{C1.2}, we have the Poincar\'e inequality
$$\mu({\hat f}^2)\le \tt C\E(f,f),\ \ f\in\D(\E)$$ for some constant $\tt C>0$. Combining this with (\ref{*0}) we prove (\ref{PH}) for
$C= C_1+\tt C(c_\phi+C_2).$
\end{proof}

\section{Extensions to the sub-Markov setting}

In this section we let $(\E,\D(\E))$ be a non-conservative Dirichlet form on $L^2(\mu)$, for which either $1\notin \D(\E)$ or $1\in\D(\E)$ but $\E(1,1)>0.$ In this case the Dirichlet form is irreducible if and only if for any $f\in\D(\E)$ with  $\E(f,f)=0$ one has  $f=0.$ Let $\|P\|_p$ be the norm in $L^p(\mu)$ for $p\ge 1.$
 Below is an extension of Corollary \ref{C1.2} to the present situation.

\beg{thm}\label{T4.1}  Let $(\E,\D(\E))$ be a non-conservative
irreducible Dirichlet form on $L^2(\mu)$. Then  the Poincar\'e inequality
\beq\label{NP}\mu(f^2)\le C\E(f,f),\ \  f\in\D(\E)\end{equation}
  holds for some $C>0$ if and only if the defective Poincar\'e inequality
\beq\label{NDP} \mu(f^2)\le C_1\E(f,f)+C_2\mu(|f|)^2,\ \ f\in\D(\E)\end{equation}
 holds for some $C_1,C_2>0.$\end{thm}

\beg{proof} Assume that (\ref{NDP}) holds. We aim  to prove (\ref{NP}) for some constant $C>0$.  According to \cite[Proposition 3.2]{W03} we need only to prove the weak Poincar\'e inequality
$$\mu(f^2)\le \aa(r)\E(f,f)+r\|f\|_\infty^2,\ \ r>0, f\in\D(\E)$$ for some function $\aa: (0,\infty)\to (0,\infty).$ This is ensured by
Lemma \ref{T4.2} below.
\end{proof}

Before introducing Lemma \ref{T4.2},    we present an application of Theorem \ref{T4.1} to  sub-Markov operators.

\beg{cor}\label{C4.2}      Let $P$  be a sub-Markov operator on $L^2(\mu)$; i.e. $P$ is a contraction linear operator on $L^2(\mu)$ with $P1\le 1$ such that $f\ge 0$ implies $Pf\ge 0.$  Let  $P^*$ be the adjoint operator of $P$.  Assume that     ${\rm Ker}(1-P^*P)=\{0\}$.
  Then   $$\|P\|_2:=\sup_{\mu(f^2)\le 1}\mu((Pf)^2)^{\ff 1 2}<1$$ if and only if
$$ \|P\|_{tail} := \lim_{R\to\infty} \sup_{\mu(f^2)\le 1} \mu\big({(|Pf|-R)^+}^2\big) <1.
 $$\end{cor}

\beg{proof} It suffices to prove $\|P\|_2<1$ from $\|P\|_{tail}<1$.  To apply Theorem \ref{T4.1}, let
  $$\E(f,g)= \mu(f (1-P^*P)g)),\ \ \D(\E)= L^2(\mu).$$   Then $(\E,\D(\E))$ is a symmetric Dirichlet form.
  Since Ker$\,(1-P^*P)=0$ and noting that the contraction of $P^*$ in $L^2(\mu)$ implies
  $$\mu((f-P^*Pf)^2) =\mu(f^2)-\mu((Pf)^2) +\mu((P^*Pf)^2) -\mu((Pf)^2)\le \E(f,f),$$ this Dirichlet form is non-conservative and irreducible.
   By $\|P\|_{tail}<1$,
  there exist $R>0$ and $\vv\in (0,1)$ such that
 $$\mu\big({(|Pf|-R)^+}^2\big)\le \vv^2,\ \ \ \mu(f^2)\le 1.$$ So, for any $f$ with
  $\mu(f^2)=1,$
$$\mu((Pf)^2)\le  \mu\big(|Pf|(|Pf|-R)^+\big)  + R\mu(|Pf|)\le\vv + R\mu(|Pf|) .$$ This implies
$$\E(f,f):= \mu(f^2)- \mu((Pf)^2)\ge 1-\vv -R\mu(|f|)\ge \ff{1-\vv } 2 -\ff{R^2}{
  2(1-\vv)}\mu(|f|)^2.$$ Therefore,
$$\mu(f^2)=1\le \ff{2}{1-\vv}\,\E(f,f) +\ff{R^2}{(1-\vv)^2}\, \mu(|f|)^2,\ \ \ f\in L^2(\mu).
  $$ Thus,
  (\ref{NDP}) holds for $C_1= \ff 2 {1-\vv}, C_2= \ff{R^2}{(1-\vv)^2}.$ By Theorem \ref{T4.1}, there exists $C>0$ such that
$$\mu(f^2)\le C\E(f,f)=C\mu( f(1-P^*P)f)= C \mu(f^2) -C \mu((Pf)^2).$$
  This implies that $\|P\|_2^2\le \ff{C-1}C <1.$ \end{proof}

Correspondingly to  Remark 1.1(a) in the ergodic setting,  we would believe  that
in the present situation  $\|P\|_\tau<1$ is not enough to imply $\|P\|_2<1$.  But it seems hard to give a proof  or a counterexample.

\beg{lem}\label{T4.2}   A non-conservative Dirichlet form
$(\E,\D(\E))$ on $L^2(\mu)$ is irreducible if and only if  there exists $\aa:
(0,\infty)\to (0,\infty)$ such that
\beq\label{WP} \mu(f^2)\le \aa(r)\E(f,f) + r \|f\|_\infty^2,\ \ \ r>0, f\in \D(\E).
\end{equation}
Consequently, for any symmetric $($sub-$)$ Markov semigroup $P_t$ on
$L^2(\mu)$, $\lim_{t\to\infty} \|P_tf\|_2= 0$ for any $f\in L^2(\mu)$   if and only if
\beq\label{UC}\lim_{t\to\infty}   \|P_t \|_{\infty\to 2}=0.\end{equation}
\end{lem}

\beg{proof} (a) Let $P_t$ be the associated semigroup
of $(\E,\D(\E))$. Then $(\E,\D(\E))$ is irreducible if and only if $\mu((P_tf)^2)\to 0$ as
$t\to\infty$ for any $f\in L^2(\mu).$ Indeed, the sufficiency is obvious since $\E(f,f)=0$ implies
$$\ff{\d}{\d t}\mu(|P_tf|^2)=2\E(P_tf,P_tf)\le 2\E(f,f)=0,$$ so that $\mu(|P_tf|^2)=\mu(f^2)$ for all $t
\ge 0$,  while the necessity holds since the irreducibility of the Dirichlet form implies that
$$\lim_{t\to \infty} \mu(|P_tf|^2) = \lim_{t\to\infty} \int_{(-\infty,0)} \e^{2\ll t} \d \|E_\ll(f)\|^2=0,\ \ f\in L^2(\mu),$$ where
$\{E_\ll\}_{\ll\le 0}$ is the spectral family of the generator such that $\d\|E_\ll(f)\|^2$ is a finite measure on $(-\infty,0)$ for $f\in L^2(\mu)$
(by the irreducibility $\{0\}$ is a null-set of the measure).

On the other hand,
by \cite[Theorem 2.1]{RW} with $\Phi(f)=   \|f\|_\infty^2$,
(\ref{WP}) holds for some $\aa$ if and only if (\ref{UC}) holds.
So, the second assertion follows from the first one.

(b)  Let $f\in \D(\E)$ with $\E(f,f)=0$. Let $\tt f =
 |f| \land 1.$
  We have
$\E(\tt f, \tt f)=0$.
So, applying (\ref{WP}) to $\tt f$  we obtain
$\mu(\tt f)\le r$ for all $r>0$. This implies
$\tt f=0$   and thus, $f=0.$

(c) Now, let $(\E,\D(\E))$ be irreducible. We claim that (\ref{WP}) holds for some function
$\aa: (0,\infty)\to (0,\infty).$ Otherwise,   there exist some  $r>0$ and a sequence $\{f_n\}\subset \D(\E)$  such that
\beq\label{2.1} 1=\mu(f_n^2)> n \E(f_n,f_n) +
r \|f_n\|_\infty^2,\ \ \ n\ge 1.
\end{equation} Since
$\E(|f_n|, |f_n|)\le \E(f_n, f_n)$, we may and do assume that $f_n\ge 0$ for all $n\ge 1.$
 Since $\{f_n\}$ is bounded  in $L^2(\mu)$,
  there exist   $f\in L^2(\mu)$ and a subsequence $\{f_{n_k}\}$
  such that $f_{n_k}$ converges weakly to $f$ in $L^2(\mu)$.

Let $P_t$ be the
(sub-) Markov semigroup and $(L,\D(L))$ the generator associated to $(\E,\D(\E))$. Then $P_t f\in \D(L)$ for any $t>0.$ By the symmetry of $P_t$ and the weak convergence of $\{f_{n_k}\}$ to $f$ in $L^2(\mu)$, we have
$$ \lim_{k\to\infty} \mu((P_t f_{n_k})g)= \lim_{k\to \infty}\mu(f_{n_k} P_t g)=
\mu(f P_t g)= \mu((P_t f)g),\ \ \ g\in L^2(\mu).$$ This implies
\beq\label{2.3}\beg{split}& \lim_{k\to \infty} \E(P_t f_{n_k}, g)
=-\lim_{k\to\infty} \mu((P_t f_{n_k})Lg)\\
&= -\mu((P_t f)Lg)=\E(P_t
f,g),\ \ \ g\in \D(L).\end{split}\end{equation} Moreover, due to (\ref{2.1})
and the symmetry of $\E$,
$$\lim_{k\to \infty} \E(P_tf_{n_k}, g)^2 \le \lim_{k\to\infty} \E(P_tf_{n_k},P_t f_{n_k})\E(g,g) \le
\lim_{k\to\infty}\E(f_{n_k},f_{n_k})\E(g,g)=0.$$
 Combining this with (\ref{2.3}) we conclude that $\E(P_tf, P_tf)=0$ for all $t>0$. Thus, by the irreducibility,
  $P_t f=0$ holds for all $t>0.$ This implies $f=0$ by the strong continuity of $P_t$ in $L^2(\mu)$.
 Since (\ref{2.1}) implies $f_n\le r^{-1/2},$ by the weak convergence of $\{f_{n_k}\}$ to
 $f=0$ in $L^2(\mu)$ and $1\in L^2(\mu)$, we obtain
 $$\lim_{k\to \infty} \mu(f_{n_k}^2)\le r^{-1/2}\lim_{k\to\infty} \mu(f_{n_k})=0.$$ This  contradicts
 the assumption that $\mu(f_n^2)=1$ for all $n\ge 1.$ Therefore, (\ref{WP}) holds for some function
  $\aa: (0,\infty)\to (0,\infty).$\end{proof}

\paragraph{Acknowledgement.} I would like to thank professor Laurent Miclo  for sending me his exciting paper \cite{M} and useful discussions,
as well as  the referee, Professor L. Gross and Dr. Shaoqin Zhang for helpful comments and corrections.

\beg{thebibliography}{99}

\bibitem{A} S. Aida, \emph{Uniformly positivity improving property,
Sobolev inequalities and spectral gap,} J. Funct. Anal.
158 (1998), 152--185.

\bibitem{C} P. Cattiaux, \emph{A pathwise approach of some classical inequalities,} Pot. Anal. 20(2004), 361--394.

\bibitem{GW} F. Gong, L. Wu, \emph{Spectral gap of positive operators and
applications,} J. Math. Pures Appl. 85(2006), 151--191.

\bibitem{Gross} L. Gross, \emph{Logarithmic Sobolev inequalities,} Amer. J. Math. 97 (1975) 1061Ð1083.

\bibitem{H} M. Hino, \emph{Exponential decay of positivity preserving
semigroups on $L^p$,} Osaka J. Math. 37 (2000), 603--624.


\bibitem{K} S. Kusuoka, \emph{Analysis on Wiener spaces II:
differential forms,} J. Funct. Anal. 103 (1992), 229--274.

\bibitem{LO} R. Lata{\l}a, K. Oleszkiewicz, \emph{Between
Sobolev and Poincar\'e,} Lecture Notes in Math. Vol. 1745,
147--168, 2000.

\bibitem{LS} G. F. Lawler, A. D. Sokal, \emph{Bounds on the $L^2$ spectrum for Markov chains and Markov processes: a generalization of Cheeger's inequality,} Trans. Amer. Math. Soc. 309(1988), 557--580.

\bibitem{24} J. R. Lee, S. Oveis Gharan, L. Trevisan, \emph{Multi-way spectral partitioning and higher-order Cheeger inequalities,} arXiv: 1111.1055v4.
STOC'12 Proceedings of the 44th Symposium on Theory of Computing, pp. 1117--1130, New York, NY, USA, 2012.

\bibitem{M} L. Miclo, \emph{On hyperboundedness and spectrum of Markov operators,} http://hal.archives-ouvertes.fr/hal-00777146.

\bibitem{RS} M. Reed, B. Simon,  \emph{Methods of Modern
Mathematical Physics IV: Analysis of Operators,} Academic Press,
New York, 1978.

 \bibitem{RS1}  M. Reed, B. Simon, \emph{Method of Modern
Mathematical Physics, I: Functional Analysis,} Academic Press,
1980.

\bibitem{RW} M. R\"ockner, F.-Y. Wang, \emph{Weak Poincar\'e inequalities
and $L^2$-convergence rates of Markov semigroups,}  J. Funct.
Anal. 185(2001), 564--603.

\bibitem{SH} B. Simon, R. H$\phi$egh-Krohn, \emph{Hypercontractive
semigroups and two dimensional self-coupled Bose fields,}
J. Funct. Anal. 9(1972), 121--180.

\bibitem{W00} F.-Y. Wang, \emph{Functional inequalities, semigroup
properties and spectrum estimates,} Infin. Dimer. Anal. Quant.
Probab. Relat. Topics, 3(2000), 263--295.

\bibitem{W02} F.-Y. Wang, \emph{Functional inequalities and spectrum
estimates: the infinite measure case,} J. Funct. Anal. 194(2002),
288--310.

\bibitem{W03} F.-Y. Wang, \emph{Functional inequalities for the decay of sub-Markov semigroups,} Pot. Anal.
18(2003), 1--23.

\bibitem{W04} F.-Y. Wang, \emph{Spectral gap for hyperbounded operators,}  Proc. Amer. Math. Soc. 132(2004), 2629--2638.

\bibitem{W05}  F.-Y. Wang, \emph{A generalization of Poincar\'e
and log-Sobolev inequalities,} Potential Analysis 22(2005), 1--15.

\bibitem{Wbook}
F.-Y. Wang,   \emph{Functional Inequalities, Markov Processes and
Spectral Theory}, Science Press, Beijing, 2005.

\bibitem{Wu} L. Wu, \emph{Uniformly integrable operators and large deviations
for Markov processes,} J. Funct. Anal. 172(2000), 301--376.

\bibitem{Wu2} L. Wu, \emph{Uniform positive improvingness, tail norm condition and spectral gap,}
www.math.kyoto-u.ac.jp/probability/symp/sa01/wu/pdf.

\end{thebibliography}

\end{document}